\documentstyle[11pt]{amsart}
\newtheorem{thm}{Theorem}
\newtheorem{lem}{Lemma}
\newtheorem{prop}{Proposition}
\newtheorem{cor}{Corollary}
\newtheorem{defn}{Definition}
\newtheorem{remark}{Remark}
\theoremstyle{definition}
\def\real{{\Bbb R}}
\def\nat{{\Bbb N}}
\def\A{{\cal A}}
\def\K{{\cal K}}
\def\L{{\cal L}}
\def\N{{\cal N}}
\def\P{{\cal P}}
\def\chix{\raise.5ex\hbox{$\chi$}}
\def\varep{\varepsilon}
\newcommand{\conv}{\operatorname{conv}}
\newcommand{\diam}{\operatorname{diam}}

\begin{document}
\baselineskip=18pt
\title[Pe{\l}czy\'nski's Property (v)]{Pe{\l}czy\'nski's Property (v) on
spaces of vector-valued functions}
\author{Narcisse Randrianantoanina}
\address{Department of Mathematics, The University of Texas at Austin,
Austin, TX 78712-1082}
\email{nrandri@@math.utexas.edu}

\maketitle

\centerline{Preliminary Version}

\begin{abstract}
Let $E$ be a separable Banach space and $\Omega$ be a compact Hausdorff
space. It is shown that the space $C(\Omega,E)$ has property (V) if and
only if $E$ does. Similar result is also given for Bochner spaces
$L^p(\mu,E)$ if $1<p<\infty$ and $\mu$ is a finite Borel measure on $\Omega$.
\end{abstract}

\section{Introduction}

Let $E$ and $F$ be Banach spaces and suppose $T:E\to F$  is a bounded linear
operator. The operator $T$ is said to be unconditionally converging if $T$
does not fix any copy of $c_0$. A Banach space $E$ is said to have Pe{\l}czy\'nski's
property (V) if every unconditionally converging operator with domain $E$
is weakly compact.
In a fundamental paper \cite{PL1}, Pe{\l}czy\'nski showed that if $\Omega$
is a compact Hausdorff space then the space $C(\Omega)$, of all continuous
scalar valued functions on $\Omega$, has property (V); and he asked
(\cite{PL1} Remark~1, p.~645; see also \cite{DU} p.~183)
if for a Banach space $E$ the abstract continuous function space
$C(\Omega,E)$ has property (V).
This question has been considered by several authors.
Perhaps the sharpest result so far is in the paper of Cembranos, Kalton,
E.~Saab and P.~Saab \cite{CKSS} where they proved that if $E$ has
property (U) and contains no copy of $\ell^1$ then $C(\Omega,E)$ has
property (V).
There are however many known examples of Banach spaces that have property (V) but fail to satisfy the above conditions.
For instance, Kisliakov in \cite{K} (see also Dealban \cite{DE} independently)
showed that the disk algebra has property (V); Bourgain did the same for ball
algebras and polydisk algebras in \cite{BO5} and $H^\infty$ in \cite{BO4}. For more information and examples of spaces with property (V), we refer to \cite{GOS2} and \cite{SS7}.

In this note, we obtain a positive answer to the above question for the
separable case; namely we prove that if $E$ is a separable Banach space
then $C(\Omega,E)$ has property (V) if and only if $E$ does.
We present also some applications of the main theorem to Banach spaces
of compact operators as well as for Bochner function spaces.

Our notation is standard Banach space terminology as may be found in the
books \cite{D1} and \cite{DU}.

\section{Definitions and some preliminary results}

\begin{defn}
Let $E$ be a Banach space.
A series $\sum\limits_{n=1}^\infty x_n$ in $E$ is said to be weakly unconditionally
Cauchy {\rm (WUC)} if for every $x^*$ in $E^*$, the series $\sum\limits_{n=1}^\infty
|x^*(x_n)|$ is convergent.
\end{defn}

There are many criteria for a series to be a WUC series (see for instance
\cite{D1}).

The following proposition was proved by Pe{\l}czy\'nski in \cite{PL1}.

\begin{prop}
For a Banach space $E$, the following assertions are equivalent:
\begin{enumerate}
\item[(i)] A subset $H\subset E^*$ is relatively weakly compact whenever
$\lim\limits_{n\to\infty} \sup\limits_{x^*\in H} |x^*(x_n)| =0$ for every {\rm WUC} series
$\sum\limits_{n=1}^\infty x_n$ in $E$;
\item[(ii)] For any Banach space $F$, every bounded operator $T:E\to F$
that is unconditionally converging is weakly compact.
\end{enumerate}
\end{prop}

\begin{defn}
A subset $H\subset E^*$ is called a (V)-subset if $\lim\limits_{n\to\infty}
\sup\limits_{x^*\in H} |x^*(x_n)|=0$ for every {\rm WUC} series
$\sum\limits_{n=1}^\infty x_n$ in $E$.
\end{defn}

So a Banach space $E$ has property (V) if (and only if) every (V)-subset of
$E^*$ is relatively weakly compact.
This leads us to study (V)-subsets of the dual of $C(\Omega,E)$ for a
given Banach space $E$ and a compact Hausdorff space $\Omega$.

Recall that the space $C(\Omega,E)^*$ is isometrically isomorphic to the
Banach space $M(\Omega,E^*)$ of all weak*-regular $E^*$-valued measures
of bounded variation defined on the $\sigma$-field $\Sigma$ of Borel subsets
of $\Omega$ and equipped with norm $\|m\| = |m|(\Omega)$, where $|m|$ is the
variation of $m$.
In this section we study different structures of subsets of $M(\Omega,E^*)$.

Let us begin by recalling some classical facts:
Fix  $\lambda$ a probability measure on $\Sigma$ and let $m\in M(\Omega,E^*)$
with $|m|\le\lambda$ and $\rho$ be a lifting of $L^\infty (\lambda)$
(see \cite{DIU} and \cite{IT}).
For $x\in E$, the scalar measure $xom$ has density
$dxom/d\lambda\  \in L^\infty (\lambda)$.
We define $\rho (m)(\omega) (x) = \rho (dxom/d\lambda)(\omega)$.
It is well known that
$$x(m(A)) = \int_A \langle \rho (m)(\omega),x\rangle \,d\lambda(\omega)$$
and
$$|m|(A) = \int_A \|\rho (m)(\omega)\|\,d\lambda(\omega)$$
for every measurable subset $A$ of $\Omega$.
Note also that $\omega\mapsto \rho (m)(\omega) (\Omega\to E^*)$ is
weak*-scalarly measurable.

The following proposition can be deduced from \cite{B1} but  we will present
a direct proof for sake of completeness.

\begin{prop}
Let $H$ be a bounded subset of $M(\Omega,E^*)$.
If $H$ is a (V)-subset then $V(H)=\{|m|,m\in H\}$ is relatively
weakly compact in $M(\Omega)$.
\end{prop}

\begin{pf}
Assume that $V(H)$ is not relatively weakly compact.
Since the space $C(\Omega)$ has property (V), there exists a WUC series
$\sum\limits_{n=1}^\infty e_n$ in $C(\Omega)$, sequence $(m_n)_n$ in $H$ and
$\varep >0$ so that $\langle e_n,|m_n|\rangle \ge \varep$ for
each $n\in \nat$.
Let $\lambda= \sum\limits_{n=1}^\infty {1\over 2^n} |m_n|$.
Since $|m_n| \le 2^n\lambda$, for a lifting $\rho$ of $L^\infty(\lambda)$,
there exists a weak*-scalarly measurable map $g_n:\Omega \to E^*$ so that:
\begin{enumerate}
\item[(a)] $\|g_n(\cdot)\| \in L^\infty$ and
$|m_n| (A) = \int_A \|g_n(\omega)\|\,d\lambda (\omega)$ for all $A\in\Sigma$;
\item[(b)] $\langle e,m_n(A)\rangle = \int_A \langle e,g_n(\omega)\rangle\,
d\lambda(\omega)$ for $e\in E$ and $A\in\Sigma$;
\item[(c)] $\rho (g_n) = g_n$.
\end{enumerate}
Now since $C(\Omega,E)$ is norming for $M(\Omega,E^*)$, there exists
$\theta_n \in C(\Omega,E)$, \ with $\|\theta_n\|=1$ and such that
$\langle \theta_n,m_n\rangle \ge \|m_n\| - \varep/2$; i.e.,
$$\int \theta_n (\omega)\,d m_n(\omega) \ge \int \|g_n(\omega)\|\,
d\lambda (\omega) - {\varep\over2}$$
or
$$\int \langle \theta_n (\omega, g_n(\omega)\rangle\,d\lambda(\omega)
\ge \int \|g_n(\omega)\|\,d\lambda (\omega) - {\varep\over2}\ .$$
Notice also that since $\|\theta_n(\omega) \|\le 1$,
$\langle \theta_n(\omega),g_n(\omega)\rangle \le \|g_n(\omega)\|$
and  we get that
\begin{align}
&\Big| \int e_n(\omega) \theta_n(\omega)\,dm_n(\omega)
- \int e_n (\omega) \,d|m_n| (\omega)\Big|      \notag\\
&\qquad
= \Big| \int e_n(\omega) \left( \langle \theta_n(\omega),g_n(\omega)
\rangle - \|g_n(\omega)\|\right) \,d\lambda(\omega)\Big|\notag\\
&\qquad
\le \int \|g_n(\omega)\|\,d\lambda(\omega)
- \int\langle \theta_n(\omega),g_n(\omega)\rangle \,d\lambda(\omega)\notag\\
&\qquad \le \int \|g_n(\omega)\|\,d\lambda(\omega)
- \biggl( \int \|g_n(\omega)\|\,d\lambda(\omega)- {\varep\over2}\biggr)
= {\varep\over2}\ .\notag
\end{align}
So for each $n\in\nat$,
$$\Big|\int e_n(\omega) \theta_n (\omega)\,dm_n\Big| > {\varep\over2}\ .$$
Fix $\psi_n = e_n(\cdot)\theta_n(\cdot)$; the function $\psi_n$ belongs to
$C(\Omega,E)$ and we claim that $\sum\limits_{n=1}^\infty \psi_n$ is a WUC series
in $C(\Omega,E)$.
For that it is enough to notice that for any finite subset $\sigma$ of $\nat$,we get
\begin{align}
\Big\| \sum_{n\in\sigma} e_n(\omega)\theta_n(\omega)\Big\|_E
&= \sup_{\|x^*\|\le 1} \Big| \sum_{n\in\sigma} e_n(\omega)
\langle \theta_n(\omega),x^*\rangle\Big| \notag\\
&\le C\ \text{ (for some constant $C$).}\notag
\end{align}
Now $\langle \psi_n,m_n\rangle\ge \varep/2$,  $\forall\ n\in\nat$.
Contradiction with  the assumption that $H$ is a (V)-subset.
\end{pf}

For the next proposition, we will use the following notation:
for a given measure $m\in M(\Omega,E^*)$ and $A\in \Sigma$, $m\chix_A$
denotes the measure $(\Sigma\to E^*)$ given by $m\chix_A (B) = m(A\cap B)$
for any $B\in\Sigma$.

\begin{prop}
Let $H$ be a (V)-subset of $M(\Omega,E^*)$ and $(A_m)_{m\in H}$
--- a collection  of measurable subsets of $\Omega$.
Then the subset $\{m\chix_{A_m},m\in H\}$ is a (V)-subset of
$M(\Omega,E^*)$.
\end{prop}

\begin{pf}
Assume that $H$ is a (V)-subset of $M(\Omega,E^*)$.
By Proposition  3, $V(H)$ is relatively weakly compact in $M(\Omega)$.
Let $\lambda$ be a control measure for $V(H)$.
Fix a sequence $(m_n\chix_{A_{m_n}})_{n\in\nat}$ in
$\{m\chix_{A_m},m\in H\}$.
We need to show that the countable subset $\{m_n\chix_{A_{m_n}},n\in\nat\}$
is a (V)-subset.
Let $\sum\limits_{n=1}^\infty f_n$ be a WUC series in $C(\Omega,E)$ with
$\sup_n \|f_n\|\le 1$.
For $\varep >0$ (fixed), there exists $\delta >0$ such that if
$A\in\Sigma$, $\lambda(A)<\delta$ then $|m|(A) <\varep/2$,\
$\forall\ m\in H$; for each $n\in\nat$, choose a compact set $C_n$ and an
open set $O_n$ such that $C_n \subset A_{m_n} \subset O_n$ and
$\lambda(O_n\setminus C_n) <\delta$.
Fix a continuous function $g_n : \Omega\to [0,1]$ with $g_n(C_n)=1$
and $g_n(\Omega\setminus O_n) =0$ and let $\phi_n = g_n f_n$.
It is not difficult to see that $\sum\limits_{n=1}^\infty \phi_n$ is a WUC
series and hence $\lim\limits_{n\to\infty} \langle m_n,\phi_n\rangle =0$.
Now we have the following estimate
\begin{align}
|\langle m_n\chix_{A_{m_n}} ,f_n\rangle|
&= \Big| \int f_n \chix_{A_{m_n}} \,dm_n\Big|\notag\\
&=\Big| \int_{A_{m_n}} f_n\,dm_n\Big|\notag\\
&\le \Big| \int_{C_n} f_n\,dm_n\Big|
+\Big| \int_{A_{m_n}\setminus C_n} f_n\, dm_n\Big|\notag\\
&\le \Big|\int_\Omega \phi_n\,dm_n\Big|
+ \Big| \int_{O_n\setminus C_n} f_n\,dm_n\Big|
+ \Big| \int_{A_{m_n}\setminus C_n} f_n\,dm_n\Big|\notag\\
&\le |\langle \phi_n,m_n\rangle| + |m_n| (O_n\setminus C_n)
+ |m_n| (A_{m_n} \setminus C_n)\notag\\
&\le |\langle \phi_n,m_n\rangle| + 2|m_n| (O_n\setminus C_n)\notag\\
&\le |\langle \phi_n,m_n\rangle| + \varep.\notag
\end{align}
This implies that $\limsup\limits_{n\to\infty} |\langle
m_n\chix_{A_{m_n}}
,f_n
\rangle| \le \varep$ and since $\varep$ is arbitrary,
we conclude that
$\lim\limits_{n\to\infty} |\langle m_n\chix_{A_{m_n}}  ,f_n\rangle|=0$.
This shows that $\{m_n\chix_{A_{m_n}},n\in\nat\}$ is a
(V)-subset.
\end{pf}

If we denote by $M^\infty (\lambda,E^*)$ the set $\{m\in M(\Omega,E^*)$;
$|m|\le \lambda\}$ then we obtain the following corollary.

\begin{cor}
Let $H$ be a (V)-subset of $M(\Omega,E^*)$ and consider $\lambda$ the
control measure of $V(H)$.
For $\varep>0$ fixed, there exists $N\in\nat$ and $H_\varep$
a (V)-subset of $M(\Omega,E^*)$ with $H_\varep\subset NM^\infty
(\lambda,E^*)$ so that $H\subseteq H_\varep + \varep B$
(where $B$ denotes the closed unit ball of $M(\Omega,E^*)$.
\end{cor}

\begin{pf}
Let $g_m :\Omega\to \real_+$ be the density of $|m|$ with respect to $\lambda$.
$$\lim_{N\to \infty}
\int_{\{\omega :g_m(\omega)>N\}}
g_m (\omega)\,d\lambda (\omega) =0\ \text{ uniformly on }\ H\ .$$
Choose $N\in\nat$ so that
$$\int_{\{\omega:g_m(\omega)>N\}} g_m(\omega) \,d\lambda(\omega)
< \varep$$
and let $A_m = \{\omega,\ g_m(\omega)\le N\}$.
It is clear that $H_\varep = \{m\chix_{A_m} ,m\in H\}$ is a subset of
$NM^\infty
(\lambda,E^*)$ and is a (V)-subset by Proposition~4.
Also each measure $m$ in $H$ satisfies $m=m\chix_{A_m} + m\chix_{A_m^c}$
with$\|m\chix_{A_m^c} \| < \varep$.
\end{pf}

Our next proposition can be viewed as a generalization of Theorem~1 of
\cite{RAN3} for sequences of weak*-scalarly measurable maps.
We denote by $(e_n)$ the unit vector basis of $c_0$, $(\Omega,\Sigma,\lambda)$ a probability space  and for any  Banach space $F$,
$F_1$ stands for the closed unit ball of $F$.

\begin{prop}
Let $Z$ be a separable subspace of a real Banach space $E$ and $(f_n)_n$
be a sequence of maps from $\Omega$ to $E^*$ that are weak*-scalarly
measurable with \par \noindent $\sup\limits_n \|f_n\|_\infty \le 1$.
Let $a,b$ be real numbers with $a<b$ then:\par
There exist a sequence $g_n\in \conv \{f_n,f_{n+1},\ldots\}$, measurable
subsets $C$ and $L$ of $\Omega$ with $\lambda (C\cup L)=1$ such that
\begin{enumerate}
\item[(i)] $\omega\in C$ and $T\in \L(c_0,Z)_1$ then either
\begin{align}
&\limsup_{n\to\infty} \langle g_n(\omega), Te_n\rangle \le b
\text{ or}\notag \\
&\liminf_{n\to\infty} \langle g_n(\omega),Te_n\rangle \ge a\ ;\notag
\end{align}
\item[(ii)] $\omega\in L$, there exists $k\in\nat$ so that for each infinite
sequence $\sigma$ of zeros and ones, there exists $T\in \L(c_0,Z)_1$ such
that for $n\ge k$,
\begin{align}
&\sigma_n=1\Longrightarrow \langle g_n(\omega), Te_n\rangle \ge b\notag\\
&\sigma_n=0\Longrightarrow \langle g_n(\omega), Te_n\rangle \le a\notag.
\end{align}
\end{enumerate}
\end{prop}

The proof is a further extension of the techniques used in \cite{T2} and
\cite{RAN3}.
We will begin by introducing some notations, some of which were already used
in \cite{T2} and \cite{RAN3}.

Let $f_n:\Omega\to E^*$ be a sequence as in the statement of the proposition.
We write $u\ll f$ (or $(u_n) \ll(f_n)$) if there exists $k\in\nat$ and
$p_1<q_1<p_2<q_2<\cdots < p_n<q_n<\cdots$ so that for  $n\ge k$,
$$u_n = \sum_{i=p_n}^{q_n} \lambda_i f_i\quad ;\ \text{with}\
\lambda_i \in [0,1]\text{ and } \sum_{i=p_n}^{q_n} \lambda_i =1\ .$$
Consider $\L(c_0,Z)_1$ the closed  unit ball of $\L(c_0,Z)$ with the strong operator
topology.
It is not difficult to see (using the fact that $Z$ is separable) that
$\L(c_0,Z)_1$ is a Polish space; in particular it has a countable basis
$(O_n)_n$.
Since $\L(c_0,Z)_1$ is a metric space, we can assume that the $O_n$'s are
open balls.

The letter $\K$ will stand for the set of all (strongly) closed subsets
of $\L(c_0,Z)_1$.

\noindent
We will say that $\omega\mapsto K(\omega)(\Omega\to \K)$ is measurable
if the set $\{\omega :K(\omega) \cap O_n\ne \emptyset\}$ is a measurable
subset of $\Omega$ for every
$n\in \nat$.

Let $h_n=\sum\limits_{i=p_n}^{a_n} \lambda_i f_i$ with $\sum\limits_{i=p_n}^{q_n} \lambda_i
=1$, $\lambda_i\ge0$ and $p_1<q_1 <p_2<q_2<\cdots$;\ $V$ be an open
subset
 of $\L(c_0,Z)_1$ and $\omega\mapsto K(\omega)$
a fixed measurable map, we set
\begin{align}
&\overline{ h_n} (\omega) = \sup_{k\ge q_n} \sup \{ \langle h_n(\omega), Te_k
\rangle; T\in V\cap K(\omega)\} \tag 1\\
&\theta (h)(\omega) = \limsup_{n\to\infty} \overline{h}_n (\omega) \tag 2
\end{align}
Notice that the definition of $\overline{ h}_n$ depends on the
representation of
$h_n$ as a  block convex combination of $f_n$'s.
Similarly we set
\begin{align}
&\widetilde{ h}_n(\omega)= \inf_{k\ge q_n} \inf \{\langle h_n(\omega),
Te_k\rangle\ ;
T\in V\cap K(\omega)\}\tag 3\\
&\varphi (h) (\omega) = \liminf_{n\to\infty} \widetilde{h}_n(\omega)\
.\tag 4
\end{align}
The proof of the following lemma is just a notational adjustment of the
proof of Lemma~2 of \cite{RAN3}.

\begin{lem}
There exists $(g_n)\ll (f_n)$ such that if $(h_n)\ll (g_n)$ then \par
$\lim\limits_{n\to\infty} \|\theta (g)-\overline{h}_n\|_1 = 0$ and
$\lim_{n\to\infty} \|\varphi (g)-\widetilde{h}_n\|_1 = 0$.
\end{lem}

\noindent {\bf Main construction:}

\noindent
Fix $a<b$ and let $\tau$ be the first uncountable ordinal.
Set $h_n^0=f_n$, we construct as in \cite{RAN3} for $\alpha < \tau$, sequences $h^\alpha =
(h_n^\alpha)_n$, measurable maps $K_\alpha :\Omega\to \K$ with the following
properties:
\begin{equation}
\text{ for } \beta <\alpha < \tau\ ,\qquad h^\alpha \ll h^\beta\ ; \tag 5
\end{equation}
For $\alpha <\tau$ and $h\ll f$ with $h_n=\sum\limits_{j=p_n}^{q_n}
\lambda_i f_i$ we define
\begin{align}
&\overline{h}_{n,\ell,\alpha} (\omega)
= \sup_{k\ge q_n} \sup \{\langle h_n(\omega),Te_k \rangle ;T\in O_\ell
\cap K_\alpha (\omega)\} \tag 6\\
&\theta_{\ell,\alpha}(h) (\omega) = \limsup_{n\to\infty}
\overline{h}_{n,\ell,\alpha}
(\omega) \notag\\
&\widetilde{h}_{n,\ell,\alpha} (\omega) = \inf_{k\ge q_n} \inf \{\langle
h_n(\omega),Te_k\rangle, T\in O_\ell \cap K_\alpha (\omega)\}\notag\\
&\varphi_{\ell,\alpha} (h)(\omega) = \liminf_{n\to\infty}
\widetilde{h}_{n,\ell,\alpha}(\omega)\ .\notag
\end{align}
Then for each $\alpha$ of the form $\beta+1$ and each $h\ll h^\alpha$,
we have
\begin{align}
&\lim_{n\to\infty} \|\theta_{\ell,\beta} (h^\alpha) -
\overline{h}_{n,\ell,\beta} \|_1 = 0\notag\\
\intertext{and}
&\lim_{n\to\infty} \|\varphi_{\ell,\beta} (h^\alpha)
- \widetilde{h}_{n,\ell,\beta}\|_1=0\ .\notag
\end{align}

If $\alpha$ is limit, we set
\begin{equation}
K_\alpha (\omega) = \bigcap_{\beta <\alpha} K_\beta (\omega)\tag 7\ .
\end{equation}
If $\alpha = \beta+1$
\begin{equation}
K_\alpha (\omega) = \{T\in K_\beta (\omega),\ T\in O_\ell \Rightarrow
\theta_{\ell,\beta} (h^\alpha) (\omega) \ge b\ ,\
\varphi_{\ell,\beta} (h^\alpha) (\omega)\le a\}\ .\tag 8
\end{equation}
The construction is done in the same manner as in \cite{RAN3} and is a
direct application of Lemma~1.

As in \cite{RAN3}, one can fix an ordinal $\alpha <\tau$ such that for a.e.
$\omega\in\Omega$,\
$K_\alpha (\omega) = K_{\alpha+1} (\omega)\ .$

Let $h= h^{\alpha+1}$,
\begin{align}
C&= \{\omega :K_\alpha (\omega) = \emptyset\}\ \text{ and} \notag\\
M&= \{\omega :K_\alpha (\omega) = K_{\alpha+1}(\omega)\ne\emptyset\}. \notag
\end{align}
Clearly $C$ and $M$ are measurable and $\lambda (C\cup L)=1$.

The next lemma is the analogue of Lemma 4 of \cite{RAN3}.

\begin{lem}
Let $\omega \in C$ and $T\in \L(c_0,Z)_1$.
If $u\ll h$ then either
\begin{align}
&\limsup_{n\to\infty} \langle u_n(\omega),Te_n\rangle \leq b\ \text{ or}\notag\\
&\liminf_{n\to\infty} \langle u_n(\omega),Te_n\rangle \geq a\notag
\end{align}
\end{lem}

\begin{pf}
Let $\omega\in C$, $T\in \L(c_0,Z)_1$ and fix $u\ll h\ll f$
(say $u = \sum\limits_{j=a_n}^{b_n} \alpha_j f_j)$;
let $S: c_0\to c_0$ be an operator defined as follows $Se_{b_n} =e_n$ and
$Se_j=0$ if $j\ne b_n$, $n\in\nat$. The operator
$S$ is obviously bounded linear with $\|S\|=1$. So $T\circ S\in \L(c_0,Z)_1=
K_0(\omega)$.
Since $T\circ S\notin K_\alpha (\omega)$, there exists a least ordinal
$\beta$ for which $T\circ S\notin K_\beta (\omega)$.
The ordinal $\beta$ cannot be a limit so $\beta = \gamma+1$ and
$T\circ S\in K_\gamma (\omega)$.
By the definition of $K_\beta(\cdot)$, there exists $\ell\in\nat$ with
$T\circ S\in O_\ell$ but either $\theta_{\ell,\gamma}(h^\beta)(\omega) \le b$
or $\varphi_{\ell,\gamma} (h^\beta)(\omega) \ge a$.
Now since $u\ll h^\beta$, we get that either
\begin{align}
&\limsup_{n\to\infty} \langle u_n(\omega), Te_n\rangle
= \limsup_{n\to\infty} \langle u_n(\omega),T\circ Se_{q_n}\rangle \notag\\
&\qquad \le \theta_{\ell,\gamma} (u)(\omega)
\le \theta_{\ell,\gamma} (h^\beta) (\omega)\le b\notag\\
\intertext{or}
&\liminf_{n\to\infty} \langle u_n(\omega), Te_n\rangle
\ge \varphi_{\ell,\gamma} (u)(\omega)
\ge \theta_{\ell,\gamma} (h^\beta)(\omega) \ge a\ .\notag
\end{align}
The lemma is proved.
\end{pf}

The following property of the measurable subset $M$ is somewhat stronger
than that obtain in Lemma~5 of \cite{RAN3} and is the main adjustment of
the entire proof.

\begin{lem}
There exists a subsequence $(n(i))$ of integers such that for a.e.
$\omega\in M$, if $\sigma$ is an infinite sequence of zeros and ones then
there exists an operator $T\in \L(c_0,Z)_1$ (which may depend on $\omega$
and $\sigma$) such that:
\begin{align}
&\sigma_i =1\Longrightarrow \langle h_{n(i)} (\omega), Te_i\rangle \ge b\notag\\
&\sigma_i=0\Longrightarrow \langle h_{n(i)}(\omega),Te_i\rangle \le a\ .\notag
\end{align}
\end{lem}

\begin{pf}
Let us denote by $F$ the set of finite sequences of zeros and ones and
$F^\infty$ the set of infinite sequences of zeros and ones.
For $s\in F$, $|s|$ will denote the length of $s$.
Let $s= (s_1,\ldots,s_n)$ and $r= (r_1,\ldots,r_m)$ with  $n\le m$. We say that
$s<r$ if  $s_i=r_i$ for $i\le m$.
Let us fix a representation of $(h_n)$ as block convex combination of
$(f_n)$:
$$h_n = \sum_{i=p_n}^{q_n} \lambda_i f_i\ .$$
We will construct sequences of integers $n(i)$ and $m(i)$; measurable
sets $B_i\subset M$ and measurable maps $Q(s,\cdot) :M\to \nat$ (for $s\in F$)
such that:
\begin{align}
&q_{n(1)} < m(1) < q_{n(2)} < m_2 <\cdots < q_{n(i)} < m(i) <\cdots \ ;
\tag 9\\
&\forall\ s\in F\ ,\ \sup \{Q(s,\omega) ; \omega\in M\} <\infty\ ;\tag {10}\\
&\lambda (M\setminus B_i)  \le 2^{-i}\ ;\tag {11}
\end{align}
\begin{align}
&\text{For all }s\in F\ ,\ \diam (O_{Q(s,\omega)}) \le {1\over |s|}\ ;\tag {12}\\
&\text{For }s,r\in F\ ,\ s<r\ \text{ and }\tag {13}\\
&\qquad \omega \in \bigcap_{|s|\le i\le |r|} B_i\ ,\ \text{ one has }\
O_{Q(r,\omega)} \subset O_{Q(s,\omega)}\ ;\notag
\end{align}
\begin{align}
&\forall\ \omega\in M\ ,\ s\in F\ ,\ K_\alpha (\omega)\cap
O_{Q(s,\omega)}\ne \emptyset\ ;\tag {14}\\
&\forall\ s\in F\ ,\ \forall \ i \le p= |s|\ ,\
\forall\ \omega \in \bigcap_{i\le j\le |s|} B_j\ ,\tag {15}\\
&s_i = 1\Longrightarrow \forall\ T\in O_{Q(s,\omega)}\ ,\
\sup_{q_{n(i)} \le k\le m(i)} \langle h_{n(i)} (\omega), Te_k\rangle
\ge b\notag\\
&s_i =0 \Longrightarrow \forall\ T\in O_{Q(s,\omega)}\ ,\
\inf_{q_{n(i)} \le k\le m(i)} \langle h_{n(i)} (\omega),Te_k \rangle
\le a\ .\notag
\end{align}

The construction is done in a similar fashion as in \cite{RAN3}; the only
difference is on the selection of the measurable map $Q(s,\cdot):\Omega\to
\nat$ so that (12) is satisfied.
For that we consider instead of $\nat$, the subset $\cal M\subset \nat$
defined by
$$\cal M = \left\{ k\in\nat \ ,\ \diam O_k \le {1\over |s|}\right\}$$
and since $\nat$ and $\cal M$ are equipped with the discrete topology, we can
replace $\nat$ by $\cal M$ and use the same argument to get
$Q(s,\cdot) :\Omega\to \cal M \cup\{0\}$.

To complete the proof, let $L= \bigcup_k \bigcap_{i\ge k} B_i$.
It is clear that $\lambda (M\setminus L)=0$. Fix $\omega\in\bigcap_{i\ge
i} B_i$
 and $\sigma \in F^\infty$, denote $\sigma = (\sigma_i)_{i\in\nat}$.
Let
$$\sigma^{(m)}= (\sigma_1,\ldots,\sigma_m) \in F\ ,\ \forall\ m\in \nat\ ;$$
By (15), we get that for $m\in\nat$ and
 $\forall\ i\le m\ ,\ \forall \ \omega\in \bigcap_{i\le j\le m} B_j $,
\begin{align*}
&\sigma_i =1 \Longrightarrow \forall\ T\in O_{Q(\sigma^{(m)},\omega)}\ ,\
\sup_{q_{n(i)} \le k\le m(i)} \langle h_{n(i)} (\omega),Te_k\rangle\ge b
\\
&\sigma_i = 0 \Longrightarrow \forall\ T\in O_{Q(\sigma^{(m)},\omega)}\
,\
\inf_{q_{n(i)}\le k\le m(i)} \langle h_{n(i)} (\omega,Te_k\rangle\le a\ .
\end{align*}
It is easy to check that the same conclusion holds for
$T\in \overline{O}_{Q(\sigma^{(m)},\omega)}$ (the closure of
$O_{Q(\sigma^{(m)}, \omega)}$ for the strong operator topology).
So if we let $\A= \bigcap_{m\in\nat}
\overline{O}_{Q(\sigma^{(m)},\omega)}$,
$\A\ne \emptyset$; in fact
$(\overline{O}_{Q(\sigma^{(m)},\omega)})_{m\in\nat}$ is a nested sequence
of
nonempty closed sets (by  (13)) of a complete metric space and such that
$\diam (\overline{O}_{Q(\sigma^{(m)},\omega)}) \to 0$ $(m\to\infty)$ (by (12))
so $\A \ne\emptyset$ (see for instance \cite{MU} p.~270).

It is now clear that if $\omega\in \bigcap_{i\ge k}B_i$ and $T\in\A$,
for $i \ge k$,
\begin{align}
&\sigma_i=1\Longrightarrow \sup_{q_{n(i)}\le k\le m(i)}
\langle h_{n(i)} (\omega),Te_k\rangle \ge b\notag\\
&\sigma_i=0\Longrightarrow \sup_{q_{n(i)}\le k\le m(i)}
\langle h_{n(i)} (\omega),Te_k\rangle \le a\ .\notag
\end{align}
We complete the proof as in \cite{RAN3}:
choose $k(i) \in [q_{n(i)} ,m(i)]$ such that
$$\sup_{q_{n(i)} \le k\le m(i)} \langle h_{n(i)} (\omega),Te_k\rangle
= \langle h_{n(i)} (\omega),Te_{k(i)}\rangle$$
for $\sigma_i=1$ and
$$\inf_{q_{n(i)}\le k\le m(i)} \langle h_{n(i)}(\omega),Te_k\rangle
= \langle h_{n(i)} (\omega),Te_{k(i)}\rangle$$
for $\sigma_i=0$.   The sequence
$(k(i))$ is an increasing sequence by (9) so one can construct an operator
$S:c_0\to c_0$ with  $Se_i = e_{k(i)}$ \ $\forall\ i\in\nat$ and it is now
clear that if
\begin{align}
&\sigma_i =1\Longrightarrow \langle h_{n(i)}(\omega),T\circ Se_i\rangle\ge b
\notag\\
&\sigma_i=0\Longrightarrow \langle h_{n(i)}(\omega),T\circ Se_i\rangle\le a\ .
\notag
\end{align}
The proof of the lemma is complete.
For the proposition, we take $g_i = h_{n(i)}$ for $i\in \nat$.
\end{pf}

\section{Main Theorem}

\begin{thm}
Let $E$ be a separable Banach space and $\Omega$ be a compact Hausdorff space.
Then the space  $C(\Omega,E)$ has property (V) if and only if $E$ has
property (V).
\end{thm}

\begin{pf}
If $C(\Omega,E)$ has property (V), then the space $E$ has property (V) since
$E$ is isomorphic to a complemented subspace of $C(\Omega,E)$.
Conversely, assume that $E$ has property (V).
Let $H$ be a (V)-subset of $M(\Omega,E^*)$.
Our goal is to show that $H$ is relatively weakly compact.
Using Corollary~1, one can assume without loss of generality that there
exists a probability measure $\lambda$ on $\Sigma$ such that $|m|\le\lambda$
for each $m\in H$.
Observe that if $E$ has property (V), then $E^*$ is weakly sequentially
complete and thus $M(\Omega,E^*)$ is weakly sequentially complete as shown
in \cite{T2} (Theorem~17).
If $H$  is not relatively weakly compact, then it contains a sequence $(m_n)_n$
that is equivalent to the $\ell^1$-basis.
By Theorem~14 of \cite{T2}, there exists $m'_n\in \conv \{m_n,m_{n+1},\ldots\}$
and $\Omega'\subset\Omega$, $\lambda(\Omega')>0$ so that for $\omega\in
\Omega'$, there exists $\ell\in\nat$ such that
$(\rho (m'_n)(\omega))_{n\ge\ell}$ is equivalent to the
$\ell^1$-basis in $E^*$.
Let
\begin{equation}
f_n = \rho (m'_n)(\omega) \chix_{\Omega'}(\omega)\ ,\qquad n\in\nat\ .\notag
\end{equation}
$(f_n)_n$ is a sequence of weak*-scalarly measurable maps and
$\sup_n \|f_n\|_\infty <\infty$.

\begin{prop}
There exist a sequence $g_n\in \conv\{f_n,f_{n+1},\ldots\}$, a positive number
$\delta$ and a strongly measurable map $T:\Omega\to \L(c_0,E)_1$ such that
$$\liminf_{n\to\infty} \Big| \int \langle g_n(\omega),T(\omega)e_n\rangle
\Big| \ge \delta.$$
\end{prop}

For the proof of the proposition, let $(a(k),b(k))_{k\in\nat}$ be an
enumeration of all pairs of rationals with $a<b$.
By induction, we construct sequences $(g^k)$, measurable sets $C_k$, $L_k$
of $\Omega$ satisfying the following:
\begin{enumerate}
\item[(i)] $g^{k+1}\ll g^k$ for each $k\in \nat$;
\item[(ii)] $C_{k+1}\subset C_k$, $L_k\subset L_{k+1}$,
$\lambda (C_k\cup L_k) =1$
\item[(iii)] $\forall\ \omega\in C_k$, $\forall\ j\ge k$,
and $T\in \L(c_0,E)_1$,then  either
\begin{align}
&\limsup_{n\to\infty} \langle g_n^j(\omega),Te_n\rangle\le b(k)\notag\\
&\liminf_{n\to\infty} \langle g_n^j(\omega),Te_n\rangle\ge a(k)\ ;\notag
\end{align}
\item[(iv)] $\forall\ \omega\in L_k$, there exists $j\in \nat$ such that for
each infinite sequence $\sigma$ of zeros and ones, there
exists $T\in\L(c_0, E)_1$ such that if $n\ge j$
\begin{align}
&\sigma_n=1\Rightarrow \langle g_n^k(\omega),Te_n\rangle \ge b(k)\
\text{or} \notag\\
&\sigma_n=0\Rightarrow \langle g_n^k(\omega),Te_n\rangle \le a(k).\
\notag
\end{align}
\end{enumerate}
This is just an application of Proposition 4 inductively starting from
$g^0=f$.

Let $\P=\{k\in\nat,\  b(k)>0\}$ and $\N = \{k\in \nat,\ a(k)<0\}$.
It is clear that $\nat = \N\cup\P$.
Consider $C=\bigcap_k C_k$ and $L= \bigcup_k L_k$; we have
$\lambda (C\cup L)=1$.

\noindent{Case 1:} $\lambda (L)>0$

Since $L= \bigcup_k L_k$, there exists $k\in\nat$ such that $\lambda (L_k)>0$.
Let $(g_n) = (g_n^k)$.
If $k\in \P$ (i.e., $b(k)>0$) fix $\sigma = (1,1,1,\ldots)$.
For each $\omega\in L_k$, there exists $T\in \L(c_0,E)_1$ such that
$\langle g_n(\omega),Te_n\rangle \ge b(k)$ \ \ $\forall\ n\ge j$.
We can choose the above operator measurably using the following lemma:

\begin{lem}
There exists a strongly measurable map $T:\Omega \to \L(c_0,E)_1$ such that:
\begin{enumerate}
\item[($\alpha$)] $T(\omega)=0$ $\forall\ \omega\notin L_k$;
\item[($\beta$)] $\omega\in L_k$, there exists $j\in\nat$ such that if
$n\ge j$, $\langle g_n(\omega),T(\omega)e_n\rangle \ge b(k)$.
\end{enumerate}
\end{lem}

To see the lemma, consider $\L(c_0,E)_1$ with the strong operator topology
and $E_1^*$ with the weak*-topology.
The space $E_1^*$ is a compact metric space and hence is a Polish space.
The space $E_1^{*\nat}\times \L(c_0,E)_1$ equipped with the product topology
is a Polish space.
Let $\A$ be the following subset of $E_1^{*\nat}\times \L(c_0,E)_1$:
$$\{(x_n^*),T\}\in \A\Leftrightarrow \forall\ j\in\nat\ ,\
\langle x_n^*,Te_n\rangle \ge b(k)\ \forall\ n\ge j\ ;$$
The set $\A$ is clearly a Borel subset of $E^{*\nat}\times \L(c_0,E)_1$ and if
$\Pi :E_1^{*\nat}\times \L(c_0,E)_1\to E_1^{*\nat}$ is the first projection,
$\Pi(\A)$ is an analytic subset of $E_1^{*\nat}$.
By Theorem~8.5.3 of \cite{CO}, there exists a universally measurable map
$\Theta :\Pi (\A)\to \L(c_0,E)_1$ such that the graph of $\Theta$ is a
subset of $\A$.
Notice that if $\omega\in L_k$, $(g_n(\omega))_{n\ge1}\in \Pi(\A)$.
We define
\begin{equation}
T(\omega) =
\begin{cases}
\Theta \left((g_n(\omega))_{n\ge1}\right) &\text{if $\omega\in L_k$}\\
0&\text{otherwise.}
\end{cases}\notag
\end{equation}
It is easy to check that $T$ satisfies all the requirements of the lemma.
The lemma is proved.

Back to the proof of the proposition, we have
$\langle g_n(\omega),T(\omega)e_n\rangle \ge b(k)$
$\forall\ \omega\in L_k$, and $n\ge j$; so
$\liminf\limits_{n\to\infty} \langle g_n(\omega),T(\omega)e_n\rangle \ge b(k)$
for $\omega\in L_k$, and by Fatou's lemma,
$$\liminf_{n\to\infty} \int\langle g_n(\omega),T(\omega)e_n\rangle\,
d\lambda(\omega) \ge b(k)\lambda (L_k)$$
so if $\delta = b(k)\lambda (L_k)>0$, the proof is complete for $k\in\P$.

Now if $k\in\N$, we consider $\sigma = (0,0,\ldots)$ and choose a strongly
measurable map $\omega\mapsto T(\omega)$ (using similar argument as in the
above lemma) with $T(\omega)=0$ for $\omega\notin L_k$ and for $\omega\in L_k$,
there exists $j\in \nat$ such that $\langle g_n(\omega),T(\omega)e_n\rangle
\le a(k)<0$ for $n\ge j$.
So we get that
$$\limsup_{n\to\infty} \langle g_n(\omega),T(\omega)e_n\rangle \le a(k)$$
for each $\omega\in L_k$ and hence
$$\limsup_{n\to\infty} \int \langle g_n(\omega),T(\omega)e_n\rangle\,
d\lambda (\omega) \le a(k)\lambda (L_k)$$
which implies that
$$\liminf_{n\to\infty} \Big|\int\langle g_n(\omega),T(\omega)e_n\rangle \,
d\lambda (\omega)\Big| \ge |a(k)| \lambda (L_k)$$
so the proof is complete if $\lambda (L)>0$.

\noindent{Case 2:} $\lambda (L)=0$

Since $\lambda (C\cup L)=1$, we have that $\lambda (\Omega\setminus C)=0$.
Choose a sequence $(g_n)$ so that $(g_n)\ll (g_n^k)$ for every $k\in\nat$.
By the definition of the $C_k$'s  and (iii) we have either
$\limsup\limits_{n\to\infty} \langle g_n(\omega),Te_n\rangle\le b(k)$
or
$\liminf\limits_{n\to\infty} \langle g_n(\omega),Te_n\rangle \ge a(k)$
\ $\forall\ k\in \nat$,
and therefore for each $\omega \in C$,\
$\lim\limits_{n\to\infty} \langle g_n(\omega),Te_n\rangle$ exists for every
$T\in \L(c_0,E)_1$(*).

But for $\omega\in\Omega'$, the sequence $(f_n(\omega))_n$ is equivalent to the
$\ell^1$-basis in $E^*$ and since $(g_n)\ll(f_n)$,
$(g_n(\omega))_n$ is also equivalent to the $\ell^1$-basis in $E^*$ and therefore
$\{g_n(\omega),n\ge1\}$ cannot be a (V)-subset of $E^*$, i.e., there exists
an operator $T\in \L(c_0,E)_1$ such that $\limsup\limits_{n\to\infty}
\langle g_n(\omega),Te_n\rangle >0$; but condition $(*)$ insures that the
limit exists so for each $\omega\in\Omega'$, there exists $T\in \L(c_0,E)$
such that $\lim\limits_{n\to\infty} \langle g_n(\omega),Te_n\rangle >0$.
We now choose the operator $T$ measurably using the same argument as in the
above lemma: i.e., there exists $T:\Omega \mapsto \L(c_0,E)_1$ strongly
measurable such  that $T(\omega)=0$ for $\omega\notin \Omega'$ and
$\lim\limits_{n\to\infty} \langle g_n(\omega),T(\omega)e_n\rangle >0\quad  \forall\
\omega\in\Omega'\ .$
Let
$\delta (\omega) = \lim\limits_{n\to\infty} \langle g_n(\omega),T(\omega)e_n
\rangle\ $ \ for $\omega \in \Omega'$ and $0$ otherwise.

The map $\omega\mapsto \delta(\omega)$ is measurable and we obtain
$$\lim_{n\to\infty} \int \langle g_n(\omega),T(\omega)e_n\rangle\,d\lambda
(\omega) = \int \delta (\omega) \,d\lambda (\omega) = \delta >0\ .$$
The proof of the proposition is complete.

To complete the proof of the theorem, fix $(g_n)\ll(f_n)$,
$T:\Omega\to \L(c_0,E)_1$ strongly measurable and $\delta>0$ as in
Proposition~5.
For each $n\in\nat$, let $G_n :\Sigma\to E^*$ be the measure in
$M(\Omega,E^*)$ defined by:
$$G_n(A) = \text{weak*-}\int_A g_n(\omega)\,d\lambda(\omega)\ .$$
It is clear that $G_n\in \conv\{m_n\chix_{\Omega'},m_{n+1}\chix_{\Omega'} ,
\ldots\}$ and we will show that $\{G_n,n\ge1\}$ is not a (V)-subset
of $M(\Omega,E^*)$ to get a contradiction by virtue of Proposition~3:
Since $\omega \mapsto T(\omega)e_n$ is  norm-measurable  for each $n\in\nat$,
one can choose (using Lusin's Theorem) a compact subset $\Omega''\subset
\Omega$ with $\lambda(\Omega\setminus\Omega'') < {\delta\over3}$ and such
that the map $\omega\to T(\omega)e_n$\ $(\Omega''\to E)$ is continuous for each
$n\in \nat$.

Let $\Lambda :C(\Omega'',E)\to C(\Omega,E)$ be an extension operator
(the existence of such operator is given by Theorem~21.1.4 of \cite{SE})
and consider $t_n= \Lambda (T(\cdot)e_n|_{\Omega''})$.
The series $\sum\limits_{n=1}^\infty t_n$ is a WUC series in $C(\Omega,E)$.
In fact the operator $S:c_0 \to C(\Omega'',E)$ given by
$Se=T(\cdot)e|_{\Omega''}$ is easily checked to be linear and bounded and
$t_n=\Lambda\circ S(e_n)$ so $\sum\limits_{n=1}^\infty t_n$ is a WUC series.

The following estimate concludes the proof.
\begin{align}
&\langle t_n,G_n\rangle
= \int\langle g_n(\omega),t_n(\omega)\rangle\,d\lambda (\omega)\notag\\
&\qquad = \int_{\Omega''} \langle g_n(\omega),T(\omega)e_n\rangle\,
d\lambda(\omega) + \int_{\Omega\setminus\Omega''}
\langle g_n(\omega) ,t_n(\omega)\rangle\,d\lambda (\omega)\notag
\end{align}
so
\begin{align}
&\langle t_n,G_n\rangle
- \int\langle g_n(\omega),T(\omega)e_n\rangle\,d\lambda (\omega)\notag\\
&\qquad = \int_{\Omega\setminus\Omega''}
\langle g_n(\omega),t_n(\omega)\rangle \,d\lambda (\omega)
- \int_{\Omega\setminus\Omega''} \langle g_n(\omega),T(\omega)e_n\rangle\,
d\lambda (\omega)\notag\\
&\Big| \langle t_n,G_n\rangle - \int \langle g_n(\omega),T(\omega)e_n\rangle
\, d\lambda (\omega)\Big| \le 2 {\delta\over3}\notag
\end{align}
which implies that
$$\Big| \int \langle g_n(\omega),T(\omega)e_n\rangle \,d\lambda (\omega)\Big|
\le 2{\delta\over3} + |\langle t_n,G_n\rangle|\ .$$
Hence
$$\liminf_{n\to\infty} |\langle t_n,G_n\rangle | \ge {\delta\over3}\ .$$
This of course shows that $\{G_n,n\ge1\}$ is not a (V)-set.
The theorem is proved.
\end{pf}

Theorem 1. above has the following consequences relative to Banach spaces of
compact operators. For what follows
if $X$ and $Y$ are Banach spaces, $K_{w^*} (X^*,Y)$ denotes the Banach
space of weak* to weakly continuous compact operators from $X^*$ to $Y$
equipped with the operator norm.
We have the following corollaries.

\begin{cor}
Let $X$ and $Y$ be Banach spaces.
If $X$ is injective and $Y$ is separable and has property (V) then
$K_{w^*}(X^*,Y)$ has property (V).
\end{cor}

\begin{pf}
The space $K_{w^*}(X^*,Y)$ is isometrically isomorphic to $K_{w^*}(Y^*,X)$
which is a complemented subspace of $K_{w^*} (Y^*,C(B_{X^*})) \approx
C(B_{X^*},Y)$ has property (V) by Theorem~1.
\end{pf}

\begin{cor}
Let $X$ be a $\L_\infty$-space and $Y$ a separable  Banach space with
property (V). The space $K(X^*,Y)$ has property (V).
\end{cor}

\begin{pf}
The space $K(X^*,Y)$ is isomorphic to $K_{w^*} (X^{***},Y)$ (see \cite{RUE})
and it is well known that $X^{**}$ is injective and so $K_{w^*}(X^{***},Y)$
has property (V) by Corollary~2.
\end{pf}

We now turn our attention to Bochner spaces.
In \cite{B2}, Bombal observed that if $E$ is a closed subspace of an order
continuous Banach lattice, then $L^p (\mu,E)$ has property (V) if
$1<p<\infty$ and $E$ has property (V).
Our next result shows that for the separable case, property (V) can be lifted
to the Bochner space $L^p(\mu,E)$.

\begin{thm}
Let $E$ be a separable Banach space and $(\Omega,\Sigma,\mu)$ be a finite
measure space.
If $1<p<\infty$,then  the space $L^p(\mu,E)$ has property (V) if and only if
$E$ does.
\end{thm}

\begin{pf}
Without loss of generality, we will assume that $\Omega$ is a compact
Hausdorff space, $\mu$ is a Borel measure and $\Sigma$ is the completion
of the field of Borel-measurable subsets of $\Omega$.
For $1< p<\infty$, let $q$ such that $\frac1p +\frac1q=1$.
It is a well known fact that the dual of $L^p(\mu,E)$ is isometrically
isomorphic to the space $M^q(\mu,E^*)$ of all vector measures
$F:\Sigma\to E^*$ with
$$\|F\|_q = \sup_\pi \biggl\{ \sum_{A\in\pi}
{\|F(A)\|^q \over \mu (A)^q} \mu(A)\biggr\}^{1/q} < \infty$$
(see for instance, \cite{DU}, p.~115).

Let $H$ be a (V)-subset  of $M^q (\mu,E^*)$ and assume that $H\subseteq
M^\infty (\mu,E^*)$.
Since $C(\Omega,E)\subset L^p(\mu,E)$ and $C(\Omega,E)$ has property (V)
by Theorem~1, $H$ is relatively weakly compact in $M(\Omega,E^*)$.
Let $(m_n)_n\subset H$ and $\rho$ a lifting of $L^\infty (\mu)$.
There exists $G_n\in \conv\{m_n,m_{n+1},\ldots\}$ and
$G\in M^\infty (\mu,E^*)$ such that
$\|\rho (G_n)(\omega) - \rho (G)(\omega)\|$ converges to zero for
$\mu$ a.e. $\omega$.
By the Lebesgue dominated convergence
$$\int \| \rho (G_n)(\omega) - \rho(G) (\omega)\|^q\, d\mu(\omega)
@>> n\to\infty >  0$$
But this is equivalent to say that $\|G_n-G\|_q\to0$ which proves
that $H$ is relatively weakly compact in $M^q (\mu,E^*)$
(see for instance \cite{U}).
The theorem is proved.
\end{pf}

\begin{remark}
As it was observed in \cite{SS7}, the property (V) cannot be lifted
from $E$ to $L^\infty (\mu,E)$.
In fact the space $E = (\Sigma \oplus \ell_1^n)_{c_0}$ has  property (v)
but $L^\infty (\mu,E)$ contains a complemented copy of $\ell^1$ hence
failing property (v).
\end{remark}

\bibliography{narciref}
\bibliographystyle{plain}

\end{document}